\documentclass[12pt,twoside]{amsart}           %使用article格式:单面，A4,11pt
\usepackage{amsfonts}
\usepackage{amssymb}
\usepackage{amsmath}
\usepackage{fancyhdr}

\theoremstyle{plain}
\newtheorem{thm}{Theorem}[section]
\newtheorem{theorem}[thm]{Theorem}

\theoremstyle{definition}

\newtheorem{definition}[thm]{Definition}
\newtheorem{claim}[thm]{Claim}

\newtheorem{defn-thm}[thm]{Definition-Theorem}

\numberwithin{equation}{section}
\catcode`\@=11
\def\opn#1#2{\def#1{\mathop{\kern0pt\fam0#2}\nolimits}}
\def\underrightarrow{\mathpalette\underrightarrow@}
\def\underrightarrow@#1#2{\vtop{\ialign{$##$\cr
 \hfil#1#2\hfil\cr\noalign{\nointerlineskip}%
 #1{-}\mkern-6mu\cleaders\hbox{$#1\mkern-2mu{-}\mkern-2mu$}\hfill
 \mkern-6mu{\to}\cr}}}

\def\underleftarrow{\mathpalette\underleftarrow@}
\def\underleftarrow@#1#2{\vtop{\ialign{$##$\cr
 \hfil#1#2\hfil\cr\noalign{\nointerlineskip}#1{\leftarrow}\mkern-6mu
 \cleaders\hbox{$#1\mkern-2mu{-}\mkern-2mu$}\hfill
 \mkern-6mu{-}\cr}}}
\let\amp@rs@nd@\relax
\newdimen\ex@
\newdimen\bigaw@
\newdimen\minaw@
\newdimen\minCDaw@
\newif\ifCD@
\def\minCDarrowwidth#1{\minCDaw@#1}

\def\@CD{\def\A##1A##2A{\llap{$\vcenter{\hbox
 {$\scriptstyle##1$}}$}\Big\uparrow\rlap{$\vcenter{\hbox{%
$\scriptstyle##2$}}$}&&}%
\def\V##1V##2V{\llap{$\vcenter{\hbox
 {$\scriptstyle##1$}}$}\Big\downarrow\rlap{$\vcenter{\hbox{%
$\scriptstyle##2$}}$}&&}%
\def\={&\hskip.5em\mathrel
 {\vbox{\hrule width\minCDaw@\vskip3\ex@\hrule width
 \minCDaw@}}\hskip.5em&}%
\def\verteq{\Big\Vert&&}%
\def\noarr{&&}%
\def\vspace##1{\noalign{\vskip##1\relax}}\relax\iffalse{%
\fi\let\amp@rs@nd@&\iffalse}\fi
 \CD@true\vcenter\bgroup\relax\iffalse{%
\fi\let\\=\cr\iffalse}\fi\tabskip\z@skip\baselineskip20\ex@
 &\hfill$\m@th##$\hfill\cr}
\def\@endCD{\cr\egroup\egroup}
\def\>#1>#2>{\amp@rs@nd@\setbox\z@\hbox{$\scriptstyle
 \;{#1}\;\;$}\setbox\@ne\hbox{$\scriptstyle\;{#2}\;\;$}\setbox\tw@
 \hbox{$#2$}\ifCD@
 \global\bigaw@\minCDaw@\else\global\bigaw@\minaw@\fi
 \ifdim\wd\z@>\bigaw@\global\bigaw@\wd\z@\fi
 \ifdim\wd\@ne>\bigaw@\global\bigaw@\wd\@ne\fi
 \ifCD@\hskip.5em\fi
 \ifdim\wd\tw@>\z@
 \mathrel{\mathop{\hbox to\bigaw@{\rightarrowfill}}\limits^{#1}_{#2}}\else
 \mathrel{\mathop{\hbox to\bigaw@{\rightarrowfill}}\limits^{#1}}\fi
 \ifCD@\hskip.5em\fi\amp@rs@nd@}
\def\<#1<#2<{\amp@rs@nd@\setbox\z@\hbox{$\scriptstyle
 \;\;{#1}\;$}\setbox\@ne\hbox{$\scriptstyle\;\;{#2}\;$}\setbox\tw@
 \hbox{$#2$}\ifCD@
 \global\bigaw@\minCDaw@\else\global\bigaw@\minaw@\fi
 \ifdim\wd\z@>\bigaw@\global\bigaw@\wd\z@\fi
 \ifdim\wd\@ne>\bigaw@\global\bigaw@\wd\@ne\fi
 \ifCD@\hskip.5em\fi
 \ifdim\wd\tw@>\z@
 \mathrel{\mathop{\hbox to\bigaw@{\leftarrowfill}}\limits^{#1}_{#2}}\else
 \mathrel{\mathop{\hbox to\bigaw@{\leftarrowfill}}\limits^{#1}}\fi
 \ifCD@\hskip.5em\fi\amp@rs@nd@}
\def\@CDS{\def\A##1A##2A{\llap{$\vcenter{\hbox
 {$\scriptstyle##1$}}$}\Big\uparrow\rlap{$\vcenter{\hbox{%
$\scriptstyle##2$}}$}&}%
\def\V##1V##2V{\llap{$\vcenter{\hbox
 {$\scriptstyle##1$}}$}\Big\downarrow\rlap{$\vcenter{\hbox{%
$\scriptstyle##2$}}$}&}%
\def\={&\hskip.5em\mathrel
 {\vbox{\hrule width\minCDaw@\vskip3\ex@\hrule width
 \minCDaw@}}\hskip.5em&}
\def\verteq{\Big\Vert&}
\def\novarr{&}
\def\noharr{&&}
\def\SE##1E##2E{\slantedarrow(0,18)(4,-3){##1}{##2}&}
\def\SW##1W##2W{\slantedarrow(24,18)(-4,-3){##1}{##2}&}
\def\NE##1E##2E{\slantedarrow(0,0)(4,3){##1}{##2}&}
\def\NW##1W##2W{\slantedarrow(24,0)(-4,3){##1}{##2}&}
\def\slantedarrow(##1)(##2)##3##4{%
\thinlines\unitlength1pt\lower 6.5pt\hbox{\begin{picture}(24,18)%
\put(##1){\vector(##2){24}}%
\put(0,8){$\scriptstyle##3$}%
\put(20,8){$\scriptstyle##4$}%
\end{picture}}}
\def\vspace##1{\noalign{\vskip##1\relax}}\relax\iffalse{%
\fi\let\amp@rs@nd@&\iffalse}\fi
 \CD@true\vcenter\bgroup\relax\iffalse{%
\fi\let\\=\cr\iffalse}\fi\tabskip\z@skip\baselineskip20\ex@
 &\hfill$\m@th##$\hfill\cr}
\def\@endCDS{\cr\egroup\egroup}
\newdimen\TriCDarrw@
\newif\ifTriV@

\def\@TriCDV{\TriV@true\def\TriCDpos@{6}\@TriCD}
\def\@TriCDA{\TriV@false\def\TriCDpos@{10}\@TriCD}
\def\@TriCD#1#2#3#4#5#6{%
\setbox0\hbox{$\ifTriV@#6\else#1\fi$} \TriCDarrw@=\wd0
\advance\TriCDarrw@ 24pt \advance\TriCDarrw@ -1em
\def\SE##1E##2E{\slantedarrow(0,18)(2,-3){##1}{##2}&}
\def\SW##1W##2W{\slantedarrow(12,18)(-2,-3){##1}{##2}&}
\def\NE##1E##2E{\slantedarrow(0,0)(2,3){##1}{##2}&}
\def\NW##1W##2W{\slantedarrow(12,0)(-2,3){##1}{##2}&}
\def\slantedarrow(##1)(##2)##3##4{\thinlines\unitlength1pt
\lower 6.5pt\hbox{\begin{picture}(12,18)%
\put(##1){\vector(##2){12}}%
\put(-4,\TriCDpos@){$\scriptstyle##3$}%
\put(12,\TriCDpos@){$\scriptstyle##4$}%
\end{picture}}}
\def\={\mathrel {\vbox{\hrule
   width\TriCDarrw@\vskip3\ex@\hrule width
   \TriCDarrw@}}}
\def\>##1>>{\setbox\z@\hbox{$\scriptstyle
 \;{##1}\;\;$}\global\bigaw@\TriCDarrw@
 \ifdim\wd\z@>\bigaw@\global\bigaw@\wd\z@\fi
 \hskip.5em
 \mathrel{\mathop{\hbox to \TriCDarrw@
{\rightarrowfill}}\limits^{##1}}
 \hskip.5em}
\def\<##1<<{\setbox\z@\hbox{$\scriptstyle
 \;{##1}\;\;$}\global\bigaw@\TriCDarrw@
 \ifdim\wd\z@>\bigaw@\global\bigaw@\wd\z@\fi
 \mathrel{\mathop{\hbox to\bigaw@{\leftarrowfill}}\limits^{##1}}
 }
 \CD@true\vcenter\bgroup\relax\iffalse{\fi\let\\=\cr\iffalse}\fi
 \tabskip\z@skip\baselineskip20\ex@
 \ifTriV@
 \halign\bgroup
 &\hfill$\m@th##$\hfill\cr
#1&\multispan3\hfill$#2$\hfill&#3\\
&#4&#5\\
&&#6\cr\egroup%
\else
 \halign\bgroup
 &\hfill$\m@th##$\hfill\cr
&&#1\\%
&#2&#3\\
#4&\multispan3\hfill$#5$\hfill&#6\cr\egroup \fi}
\def\@endTriCD{\egroup}

\newcounter{Myenumi}
{\begin{list}{}{\usecounter{Myenumi}%
\settowidth{\leftmargin}{2.n}\settowidth{\labelwidth}{2.n}%
\setlength{\labelsep}{0pt}}}{\end{list}}
\newcounter{Myenumii}
{\begin{list}{}{\usecounter{Myenumii}%
\settowidth{\leftmargin}{a)n}\settowidth{\labelwidth}{a)n}%
\setlength{\labelsep}{0pt}}}{\end{list}}
\newcounter{Myenumiii}
{\begin{list}{}{\usecounter{Myenumiii}%
\settowidth{\leftmargin}{iv.n}\settowidth{\labelwidth}{iv.n}%
\setlength{\labelsep}{0pt}}}{\end{list}}

{\end{list}}

{\begin{list}{}{%
\settowidth{\leftmargin}{2.n}\settowidth{\labelwidth}{2.n}%
\setlength{\labelsep}{0pt}}}{\end{list}}

\newsymbol\onto 1310

\begin{document}

% % % % % % % % % % % % % % % % % % % % % % %
%                                           %
%                 Title                     %
%                                           %
% % % % % % % % % % % % % % % % % % % % % % %

\title{Regularity criteria in weak spaces for Navier-Stokes equations in $\mathbb{R}^3$}

\author{Zhihui Cai and Jian Zhai \\ \tiny{Department of Mahthematics, Zhejiang University, Hangzhou
310027, People's Republic of China} } \maketitle

\footnotetext{This work is supported by NSFC No.10571157.}

\begin{center}
\email{czh@cms.zju.edu.cn}, \email{jzhai@zju.edu.cn}
\end{center}

%\tableofcontents                                                          %列出章节目录

%\newpage                                                                      %换新页

% % % % % % % % % % % % % % % % % % % % % % %
%                                           %
%                 Abstract                  %
%                                           %
% % % % % % % % % % % % % % % % % % % % % % %

\begin{abstract}
In this paper we establish a Serrin type regularity criterion on the
gradient of pressure in weak spaces for the Leray-Hopf weak
solutions of the Navier-Stocks equations in $\mathbb{R}^3$. It
partly extends the results of Zhou\cite{zhou} to $L_w^\gamma$ spaces
instead of $L^\gamma$ spaces.

\end{abstract}

% % % % % % % % % % % % % % % % % % % % % % %
%                                           %
%                 Introduction              %
%                                           %
% % % % % % % % % % % % % % % % % % % % % % %

\section{Introduction}
In this paper, we consider the following Cauchy problem for the
incompressible Navier-Stokes equations :
\begin{equation} \label{eq:1}
\left\{ \begin{aligned}
         \frac{\partial{u}}{\partial{t}}+u\cdot \nabla u + \nabla p = \triangle u, \\
         div\ u=0,\\
         u(x,0) = u_0(x),
                          \end{aligned} \right.
                          \end{equation}
in $\mathbb{R}^3\times (0,T)$. Here $u=u(x,t)\in \mathbb{R}^3$ is
the velocity field, $p(x,t)$ is a scalar pressure field of an
incompressible fluid at the point $(x,t)$, and $u_0(x)$ with $div\
u_0=0$ in the sense of
distribution is the initial velocity field.\\

% % % % % % % % % % % % % % % % % % % % % % %
%                                           %
%                 background                %
%                                           %
% % % % % % % % % % % % % % % % % % % % % % %

The global existence of weak solutions in time was proved by
Leray\cite{Leray} and Hopf\cite{Hopf}. However, the answer to the
problem of global regularity for the three dimensional
incompressible Navier-Stokes equations is not known. By weak
solutions of the Navier-Stokes equations, we mean the usual
Leray-Hopf solutions:
\begin{definition}
A vector field $u=u(x,t)$ on $\mathbb{R}^3\times(0,T)$ is called a
weak solution of (\ref{eq:1}) in $\mathbb{R}^3\times(0,T)$, provided
that
\begin{eqnarray*}
(a)&& \qquad u\in L^\infty(0,T;L^2(\mathbb{R}^3))\cap
L^2(0,T;H^{1,2}(\mathbb{R}^3)),\\
(b)&& \qquad div\ u=0\qquad in \quad \mathbb{R}^3\times(0,T)\\
(c)&& \qquad \int^T_0\int_{\mathbb{R}^3} \{-u\cdot \Phi_t+\nabla
u\cdot \nabla\Phi+ (u\cdot\nabla u)\cdot\Phi\}dxdt = 0
\end{eqnarray*}
for all $\Phi\in C^\infty_c(\mathbb{R}^3\times(0,T))$ with $div\
\Phi=0$ in $\mathbb{R}^3\times(0,T)$.
\end{definition}
Serrin\cite{Serrin} proved that if $u\in
L^\alpha(0,T;L^\gamma(\mathbb{R}^3))$ is a Leray-Hopf weak solution
with $2/\alpha+3/\gamma<1, 3<\gamma<\infty$, then the solution
$u(x,t)\in C^\infty(\mathbb{R}^3\times(0,T])$. And then
Sohr\cite{Sohr} solved the limit case $2/\alpha+3/\gamma=1$.\\

There are some regularity criteria in terms of $p$ or $\nabla{p}$
for the whole space. In 2004, Zhou\cite{zhou} established a final
regularity criterion in terms of $\nabla{p}$ :
$$
\nabla{p}\in L^\alpha(0,T;L^\gamma(\mathbb{R}^3))\quad with\quad
\frac{2}{\alpha}+\frac{3}{\gamma}\leq 3,\
\frac{2}{3}<\alpha<\infty,\ 1<\gamma<\infty
$$
or $\nabla{p}\in L^{2/3}(0,T;L^\infty(\mathbb{R}^3))$, or else
$\|\nabla{p}\|_{L^\infty(0,T;L^\infty(\mathbb{R}^3))}$ is
sufficiently small.\\

Also in 2004, Kim and Kozono\cite{Kim} obtained an interior
regularity criteria in weak spaces in $\Omega\times(0,T)$ under the
assumption that $\|u\|_{L^s_w(0,T;L^r_w(\Omega))}$ is sufficiently
small for some $(r,s)$ with $\frac{2}{s}+\frac{3}{r} = 1$ and $3\leq
r < \infty$. They extended the criteria of Serrin\cite{Serrin} to
the weak space-time spaces.\\

The main purpose of this paper is to establish a regularity
criterion in weak spaces instead of $L^\gamma$ in Zhou\cite{zhou}.
Our main result is

% % % % % % % % % % % % % % % % % % % % % % %
%                                           %
%                 Main Theorem              %
%                                           %
% % % % % % % % % % % % % % % % % % % % % % %

\begin{theorem} \label{thm:1.1}
Let $u_0(x)\in L^2(\mathbb{R}^3)\cap L^q(\mathbb{R}^3)$, for $q\geq
4$, and let $div\ u_0=0$ in the sense of distribution. Suppose that
$u(x,t)$ is Leray-Hopf weak solution of (\ref{eq:1}). If
$$
\nabla p\in L^\alpha(0,T;L_w^\gamma(\mathbb{R}^3))\ with\
\frac{2}{\alpha}+\frac{3}{\gamma}< 3,\ \frac{2}{3}<\alpha<\infty,
1<\gamma<\infty,
$$
or $\nabla p\in L^\alpha(0,T;L_w^\gamma(\mathbb{R}^3))$ is
sufficiently small when$\ \frac{2}{\alpha}+\frac{3}{\gamma}= 3,\
\frac{2}{3}<\alpha<\infty, 1<\gamma<\infty$
% or $\nabla p \in L^{2/3,\infty}$, or else $\|\nabla
% p\|_{L^{\infty,1}}$
, or else $\|\nabla
p\|_{L^{\frac{2}{3}}(0,T;L^\infty_w(\mathbb{R}^3)}$ is sufficiently
small, then $u(x,t)$ is a regular solution on $[0,T]$.
\end{theorem}\vspace{0.5cm}

% % % % % % % % % % % % % % % % % % % % % % %
%                                           %
%                 Remark                    %
%                                           %
% % % % % % % % % % % % % % % % % % % % % % %

Here $L^r_w(\mathbb{R}^3)$ denotes the weak
$L^r(\mathbb{R}^3)$-space
$$
L^r_w(\mathbb{R}^3) = \{v\in L^1_{loc}(\mathbb{R}^3) :
\|v\|_{L^r_w(\mathbb{R}^3)} = \sup_{\sigma>0} \sigma|\{x\in
\mathbb{R}^3 : |v(x)| > \sigma\}|^{\frac{1}{r}}<\infty\}
$$

The Lorentz space $L^{p,r}$ is defined as follows. We have $f\in
L^{p,r}$, $1\leq p \leq \infty$, if and only if
\begin{eqnarray*}
&&\|f\|_{L^{p,r}} = (\int^\infty_0 (t^{1/p}f^*(t))^r
dt/t)^{1/r}<\infty \quad when \quad 1\leq r < \infty \\
&&\|f\|_{L^{p,\infty}} = \sup_t t^{1/p}f^*(t) < \infty \quad when
\quad r=\infty
\end{eqnarray*}
where $$f^*(t)=\inf\{\sigma:m(\sigma,f)\leq t\},\
m(\sigma,f)=\mu(\{x:|f(x)|>\sigma\})$$

We have the following properties which are useful in this paper,
with equality of norms,
$$
L^{p,p}=L^p,\quad\quad\quad L^{p,\infty}=L^p_w,\quad\quad\quad
when\quad 1\leq p\leq \infty.
$$
and $$ L^{p,r_1}\subset L^{p,r_2},\quad if\quad r_1\leq r_2
$$
In particular,
$$
L^p\subset L^p_w.
$$

% % % % % % % % % % % % % % % % % % % % % % %
%                                           %
%                 Proof of Theorem          %
%                                           %
% % % % % % % % % % % % % % % % % % % % % % %

\section{Proof of theorem \ref{thm:1.1}}
Taking $\nabla div$ on both side of (\ref{eq:1}) for smooth $(u,p)$,
one can obtain
$$
-\triangle(\nabla{p})=\sum\limits_{i,j=1}^{3}\partial_i\partial_j(\nabla(u_iu_j))
$$
Therefore the Calderon-Zygmund inequality
\begin{equation}\label{eq:2}
\|\nabla{p}\|_{L^q}\leq C_1\|\ |u||\nabla{u}|\ \|_{L^q}
\end{equation}
holds for any $1<q<\infty$.\\

Multiply both side of equation (\ref{eq:1}) by $4u|u|^2$, and
integrate over $\mathbb{R}^3$:
\begin{eqnarray}
&&\frac{d}{dt}\|u\|^4_{L^4}+4\||\nabla{u}||u|\|^2_{L^2}+2\|\nabla|u|^2\|^2_{L^2}
\leq 4\int_{\mathbb{R}^3} |\nabla{p}||u|^3 dx\label{res1}\\
&&\int_{\mathbb{R}^3} |\nabla{p}||u|^3 dx =
\int_{\mathbb{R}^3}|\nabla{p}|^{1/2}|\nabla{p}|^{1/2}|u|^3 dx\\
&&(use\ H\ddot{o}lder's\ inequality\ with\ p=4,q=\frac{4}{3})\nonumber\\
&&\leq \||\nabla{p}|^{1/2}\|_{L^4}\||\nabla{p}|^{1/2}|u|^3\|_{L^{\frac{4}{3}}}\nonumber\\
&&=(\int_{\mathbb{R}^3}|\nabla{p}|^{2})^{\frac{1}{4}}
(\int_{\mathbb{R}^3}|\nabla{p}|^{\frac{2}{3}}|u|^4)^{\frac{3}{4}}
=\|\nabla{p}\|_{L^2}^{\frac{1}{2}}(\int_{\mathbb{R}^3}|\nabla{p}|^{\frac{2}{3}}|u|^4)^{\frac{3}{4}}\nonumber\\
&&(use\ inequality\ \int |fg|\leq
\|f\|_{L^{p,\infty}}\|g\|_{L^{p',1}}, 1\leq p \leq \infty)\nonumber\\
&&\leq
\|\nabla{p}\|_{L^2}^{\frac{1}{2}}\||\nabla{p}|^{\frac{2}{3}}\|_{L^{p,\infty}}^{\frac{3}{4}}
\||u|^4\|_{L^{p',1}}^{\frac{3}{4}}\ \ \ (p=\frac{3}{2}\gamma,\
p'=\frac{3\gamma}{3\gamma-2})\nonumber
\end{eqnarray}
By definition \begin{eqnarray*} &&\|f\|_{L^{p,\infty}}=\sup_t
t^{1/p}f^*(t),\ \|f\|_{L^{p,1}}=\int^\infty_0
t^{1/p}f^*(t)\frac{dt}{t}
\\
&& f^*(t)=\inf\{\sigma:m(\sigma,f)\leq t\},\
m(\sigma,f)=\mu(\{x:|f(x)|>\sigma\})
\end{eqnarray*}
Let $k>0$, then
\begin{eqnarray*}
m(\sigma,f^k)=\mu(\{x:|f^k(x)|>\sigma\})=\mu(\{x:|f(x)|>\sigma^{1/k}\})=m(\sigma^{1/k},f)\\
\Rightarrow (f^k)^*(t)=\inf\{\sigma:m(\sigma,f^k)\leq
t\}=\inf\{\sigma:m(\sigma^{1/k},f)\leq t\}=(f^*(t))^k
\end{eqnarray*}
Then
\begin{eqnarray*}
\||\nabla{p}|^{\frac{2}{3}}\|_{L^{p,\infty}}&=&\sup_t
t^{1/p}(|\nabla{p}|^{\frac{2}{3}})^*(t)\\
&=&\sup_t t^{1/p}(|\nabla{p}|^*(t))^{\frac{2}{3}}\\
&=&(\sup_t t^{\frac{3}{2p}}|\nabla{p}|^*(t))^{\frac{2}{3}}\\
&=&\|\nabla{p}\|_{L^{\frac{2p}{3},\infty}}^{\frac{2}{3}}
\end{eqnarray*}

We obtain
\begin{eqnarray}
4\int_{\mathbb{R}^3} |\nabla{p}||u|^3 dx \leq
4\|\nabla{p}\|_{L^2}^{\frac{1}{2}}\|\nabla{p}\|_{L^{\gamma,\infty}}^{\frac{1}{2}}
\||u|^4\|_{L^{p',1}}^{\frac{3}{4}}\nonumber\\
(use\ Cauchy's\ inequality\ with\ \epsilon)\nonumber\label{eq:3}\\
\leq
\epsilon\|\nabla{p}\|_{L^2}^{2}+C(\epsilon)\|\nabla{p}\|_{L^{\gamma,\infty}}^{\frac{2}{3}}
\||u|^4\|_{L^{p',1}}
\end{eqnarray}
where $4<4p'=\frac{12\gamma}{3\gamma-2}<12$ as
$1<\gamma<\infty$.\\

Next, we want to estimate $\||u|^4\|_{L^{p',1}}$.

\begin{claim} \label{prop:1}
$$
\|u\|_{L^{\frac{12\gamma}{3\gamma-2}}}^4\leq
\|u\|_{L^4}^{4(1-\frac{1}{\gamma})}\|u\|_{L^{12}}^{\frac{4}{\gamma}}
\ ,\ 1<\gamma<\infty.$$
\end{claim}
In fact,
$u^{\frac{12\gamma}{3\gamma-2}}=u^{\frac{12(\gamma-1)}{3\gamma-2}}u^{\frac{12}{3\gamma-2}}$,
using H$\ddot{o}$lder inequality with
$p=\frac{3\gamma-2}{3(\gamma-1)},\ q=3\gamma-2$, one can easily
obtain this inequality .

\begin{claim} \label{prop:2}
$\forall\ f\in L^{p,1}(\mathbb{R}^3)$, $\forall$ $0<p_1<p<p_2$,
$$\|f\|_{L^{p,1}}\leq C(p,p_1)\|f\|_{L^{p_1}} +
C(p,p_2)\|f\|_{L^{p_2}}$$
\end{claim}
Proof.
\begin{eqnarray*}
    \|f\|_{L^{p,1}}&=&\int_0^\infty t^{1/p} f^*(t)\frac{dt}{t}=\int_0^1 t^{1/p}f^*(t)\frac{dt}{t} + \int_1^\infty t^{1/p}f^*(t)\frac{dt}{t} \\
                   &\leq&
                   \sup_t(t^{1/p_2}f^*(t))\int_0^1t^{\frac{1}{p}-\frac{1}{p_2}-1}dt+\sup_t(t^{1/p_1}f^*(t))\int_0^1t^{\frac{1}{p_1}-\frac{1}{p}-1}dt
                   \\
    &=&\frac{pp_1}{p-p_1}\|f\|_{L^{p_1,\infty}}+\frac{pp_2}{p_2-p}\|f\|_{L^{p_2,\infty}}
\end{eqnarray*}
We note that $L_p\subset L_p^*$, this completes the proof.\\

Let $q(\gamma)=\frac{12\gamma}{3\gamma-2}$,
$0<\gamma_1<\gamma<\gamma_2$. As $q(\gamma)$ is a decreasing
function, $q(\gamma_1)>q(\gamma)>q(\gamma_2)$. Due to Claim
\ref{prop:2}
\begin{eqnarray}\label{eq:4}
\||u|^4\|_{L^{p',1}}&=&\||u|^4\|_{L^{\frac{3\gamma}{3\gamma-2},1}}\leq
C(\gamma,\gamma_1)\||u|^4\|_{L^{\frac{3\gamma_1}{3\gamma_1-2}}}+C(\gamma,\gamma_2)\||u|^4\|_{L^{\frac{3\gamma_2}{3\gamma_2-2}}}\nonumber\\
&=&C(\gamma,\gamma_1)\|u\|_{L^{\frac{12\gamma_1}{3\gamma_1-2}}}^4+C(\gamma,\gamma_2)\|u\|_{L^{\frac{12\gamma_2}{3\gamma_2-2}}}^4
% &\leq&C(\gamma,\gamma_1)\|u\|_{L^{\frac{12\gamma_1}{3\gamma_1-2}}}^4+C(\gamma,\gamma_2)\|u\|_{L^4}^{4(1-\frac{1}{\gamma_2})}\|u\|_{L^{12}}^{\frac{4}{\gamma_2}}
\end{eqnarray}

% % % % % % % % % % % % % % % % % % % % % % %
%                                           %
%                 Case 1                    %
%                                           %
% % % % % % % % % % % % % % % % % % % % % % %

\textbf{Case 1:} $\frac{2}{\alpha}+\frac{3}{\gamma}= 3$($\nabla p$
is sufficiently small in $L^\alpha(0,T;L_w^\gamma(\mathbb{R}^3))$).\\

Let $\gamma_1=1$, then $q(\gamma_1)=12$. Use Claim \ref{prop:1}, We
obtain
\begin{eqnarray*}
&&\||u|^4\|_{L^{p',1}}\leq
C(\gamma,1)\|u\|_{L^{12}}^4+C(\gamma,\gamma_2)\|u\|_{L^{\frac{12\gamma_2}{3\gamma_2-2}}}^4\\
&&\leq C(\gamma,1)\|u\|_{L^{12}}^4+C(\gamma,\gamma_2)\|u\|_{L^4}^{4(1-\frac{1}{\gamma_2})}\|u\|_{L^{12}}^{\frac{4}{\gamma_2}}\\
\end{eqnarray*}
Since
\begin{eqnarray*}
\|u\|^4_{L^{12}}=\||u|^2\|^2_{L^6}\leq C\||\nabla{u}||u|\|^2_{L^2}
\end{eqnarray*}

Apply (\ref{eq:2}) and Claim \ref{prop:1} to (\ref{eq:3}), we have

% and
% $\frac{2\gamma_2}{3(\gamma_2-1)}<\frac{2\gamma}{3(\gamma-1)}\leq\alpha$,

\begin{eqnarray*}
&&\frac{d}{dt}\|u\|^4_{L^4}+4\||\nabla{u}||u|\|^2_{L^2}+2\|\nabla|u|^2\|^2_{L^2}
\leq 4\int_{\mathbb{R}^3} |\nabla{p}||u|^3 dx\\
&&\leq\epsilon
C\||\nabla{u}||u|\|^2_{L^2}+C(\epsilon)\|\nabla{p}\|_{L^{\gamma,\infty}}^{\frac{2}{3}}
(C(\gamma,1)\|u\|_{L^{12}}^4+C(\gamma,\gamma_2)\|u\|_{L^4}^{4(1-\frac{1}{\gamma_2})}\|u\|_{L^{12}}^{\frac{4}{\gamma_2}})\\
&&\leq\epsilon C\||\nabla{u}||u|\|^2_{L^2} +
C_1(\epsilon,\gamma)\|\nabla{p}\|_{L^{\gamma,\infty}}^{\frac{2}{3}}\||\nabla{u}||u|\|^2_{L^2}
+C_2(\epsilon,\gamma)\|\nabla{p}\|_{L^{\gamma,\infty}}^{\frac{2}{3}}\|\|u\|_{L^4}^{4(1-\frac{1}{\gamma_2})}\|u\|_{L^{12}}^{\frac{4}{\gamma_2}}\\
&&(use\ Cauchy's\ inequality\ with\ \delta)\\
&&\leq\epsilon C\||\nabla{u}||u|\|^2_{L^2} +
C_1(\epsilon,\gamma)\|\nabla{p}\|_{L^{\gamma,\infty}}^{\frac{2}{3}}\||\nabla{u}||u|\|^2_{L^2}+
C_3(\epsilon,\gamma,\delta)\|\nabla{p}\|_{L^{\gamma,\infty}}^{\frac{2\gamma_2}{3(\gamma_2-1)}}\|u\|^4_{L^4}+\delta\|u\|^4_{L^{12}}
\end{eqnarray*}
After choosing suitable $\epsilon$ and $\delta$, we have
\begin{equation}\label{main1}
\frac{d}{dt}\|u\|^4_{L^4}\leq
C\|\nabla{p}\|_{L^{\gamma,\infty}}^{\frac{2\gamma_2}{3(\gamma_2-1)}}\|u\|^4_{L^4}
\end{equation}

note that
$\frac{2\gamma_2}{3(\gamma_2-1)}<\frac{2\gamma}{3(\gamma-1)}<\alpha$.
Due to the integrability of $\nabla{p}$, it follows that
\begin{equation}\label{estimate}
\sup_{0\leq t\leq T}\|u(\cdot,t)\|^4_{L^4}\leq C(T)\|u_0\|^4_{L^4}.
\end{equation}

% % % % % % % % % % % % % % % % % % % % % % %
%                                           %
%                 Case 2                    %
%                                           %
% % % % % % % % % % % % % % % % % % % % % % %

\textbf{Case 2:} $\frac{2}{\alpha}+\frac{3}{\gamma}< 3$ ($\nabla p$
is bounded in $L^\alpha(0,T;L_w^\gamma(\mathbb{R}^3))$).\\

Let $\gamma_1=\frac{3\alpha}{3\alpha-2}$, then
$\frac{2}{\alpha}+\frac{3}{\gamma_1}= 3$. It follows from
(\ref{eq:4}) that
\begin{eqnarray*}
\||u|^4\|_{L^{p',1}}\leq
C(\gamma,\gamma_1)\|u\|_{L^{\frac{12\gamma_1}{3\gamma_1-2}}}^4+C(\gamma,\gamma_2)\|u\|_{L^{\frac{12\gamma_2}{3\gamma_2-2}}}^4
\end{eqnarray*}
The same as \textbf{Case 1}, we can get
\begin{eqnarray*}
&&\frac{d}{dt}\|u\|^4_{L^4}+4\||\nabla{u}||u|\|^2_{L^2}+2\|\nabla|u|^2\|^2_{L^2}
\leq 4\int_{\mathbb{R}^3} |\nabla{p}||u|^3 dx\\
&&\leq\epsilon
C\||\nabla{u}||u|\|^2_{L^2}+C(\epsilon)\|\nabla{p}\|_{L^{\gamma,\infty}}^{\frac{2}{3}}
(C(\gamma,\gamma_1)\|u\|_{L^{\frac{12\gamma_1}{3\gamma_1-2}}}^4+C(\gamma,\gamma_2)\|u\|_{L^{\frac{12\gamma_2}{3\gamma_2-2}}}^4)\\
&&\leq \epsilon C\||\nabla{u}||u|\|^2_{L^2} +(
C_2(\epsilon,\gamma,\delta_1)\|\nabla{p}\|_{L^{\gamma,\infty}}^{\frac{2\gamma_1}{3(\gamma_1-1)}}
+C_2(\epsilon,\gamma,\delta_2)\|\nabla{p}\|_{L^{\gamma,\infty}}^{\frac{2\gamma_2}{3(\gamma_2-1)}}
)\|u\|^4_{L^4} + (\delta_1+\delta_2)\|u\|^4_{L^{12}}
\end{eqnarray*}
Then, with suitable $\epsilon$ , $\delta_1$ and $\delta_2$, we have
\begin{equation}\label{main2}
\frac{d}{dt}\|u\|^4_{L^4}\leq
(C_1\|\nabla{p}\|_{L^{\gamma,\infty}}^{\alpha} +
C_2\|\nabla{p}\|_{L^{\gamma,\infty}}^{\frac{2\gamma_2}{3(\gamma_2-1)}})\|u\|^4_{L^4}
\end{equation}
As above, we can get (\ref{estimate}).\\

% % % % % % % % % % % % % % % % % % % % % % %
%                                           %
%                 Case 3                    %
%                                           %
% % % % % % % % % % % % % % % % % % % % % % %
\textbf{Case 3:} $(\alpha,\gamma)=(\frac{2}{3},\infty)$ ($\|\nabla
p\|_{L^{\frac{2}{3}}(0,T;L^\infty_w(\mathbb{R}^3)}$ is sufficiently
small).\\
Taking the limit case in (\ref{eq:3}), we have
\begin{eqnarray*}
4\int_{\mathbb{R}^3} |\nabla{p}||u|^3 dx \leq
4\|\nabla{p}\|_{L^2}^{\frac{1}{2}}\|\nabla{p}\|_{L^{\infty,\infty}}^{\frac{1}{2}}
\||u|^4\|_{L^{1,1}}^{\frac{3}{4}}\\
\leq
\epsilon\|\nabla{p}\|_{L^2}^{2}+C(\epsilon)\|\nabla{p}\|_{L^{\infty}}^{\frac{2}{3}}
\||u|^4\|_{L^{1}}\\
\leq \epsilon C \||\nabla{u}||u|\|^2_{L^2} +
C(\epsilon)\|\nabla{p}\|_{L^{\infty}}^{\frac{2}{3}}\|u\|_{L^4}^4
\end{eqnarray*}
The same as \textbf{Case 1}, after choosing suitable $\epsilon$, then use Gronwall inequality, we can get (\ref{estimate}).\\

This apriori estimate (\ref{estimate}) is what we want. Then we use
a result of Giga \cite{Giga}:
\begin{theorem}\label{Giga}
Suppose $u_0\in L^s(\mathbb{R}^3)$, $s\geq 3$. Then there exists
$T_0$ and a unique classical solution $u\in BC([0,T_0);
L^s(\mathbb{R}^3))$. Moreover, let $(0, T_*)$ be the maximal
interval such that u solves (\ref{eq:1}) in $C((0,T_*);
L^s(\mathbb{R}^3))$, $s>3$. Then
\begin{equation}\label{giga2}
\|u(\cdot,\tau)\|_{L^s}\geq \frac{C}{(T_*-\tau)^{(s-3)/2s}}
\end{equation}
with constant C independent of $T_*$ and s.
\end{theorem}

% % % % % % % % % % % % % % % % % % % % % % %
%                                           %
%         Proof of Theorem 1.1              %
%                                           %
% % % % % % % % % % % % % % % % % % % % % % %
\emph{Proof of Theorem \ref{thm:1.1}}. We follow the method of Zhou
\cite{zhou}: Since $u_0(x)\in L^2(\mathbb{R}^3)\cap
L^q(\mathbb{R}^3)$ for $q\geq 4$, due to Theorem \ref{Giga} ($s=4$),
there exists a unique solution \~{u}$(x,t)\in BC([0,T_*);
L^4(\mathbb{R}^3))$. Since $u$ is a Leray-Hopf weak solution, we
have by the uniqueness criterion of Serrin-Masuda \cite{Serrin}
\cite{Masuda}
$$
u\equiv \tilde{u} \quad on \quad [0,T_*).
$$
By the apriori estimate (\ref{estimate}), and the standard
continuation argument, we can continue our local smooth solution
corresponding to $u_0\in L^4(\mathbb{R}^3)$ to obtain $u\in
BC([0,T]; L^4(\mathbb{R}^3))\cap C^{\infty}(\mathbb{R}^3\times (0,
T])$. This completes the proof of Theorem \ref{thm:1.1}.

\vspace{1cm}
% % % % % % % % % % % % % % % % % % % % % % %
%                                           %
%                 Remark                    %
%                                           %
% % % % % % % % % % % % % % % % % % % % % % %
\emph{Remark} The limit case $(\alpha, \gamma)=(\infty, 1)$ is not
solved in this paper.

\vspace{2cm}
% % % % % % % % % % % % % % % % % % % % % % %
%                                           %
%                 Bibliography              %
%                                           %
% % % % % % % % % % % % % % % % % % % % % % %

\end{document}